\RequirePackage{fix-cm}
\RequirePackage{amsmath}
\documentclass[smallextended]{svjour3}       
\smartqed  
\usepackage{graphicx}
\usepackage[letterpaper,left=1.5in,right=1.5in,top=1in,bottom=1in]{geometry}
\usepackage[doublespacing]{setspace}
\usepackage{amsfonts} 
\usepackage{upgreek}
\usepackage{hyperref}
\usepackage[justification=centering]{caption}
\journalname{}
\begin{document}

\title{Discrete Stochastic Optimization for Public Health Interventions with Constraints
}
\subtitle{}


\author{Zewei Li  \and
       James C. Spall }


  \institute{Zewei Li \at
             Northwestern University \\
              Evanston, U.S.A\\
              \email{zeweili2025@u.northwestern.edu}           
            \and
            James C. Spall \at
            The Johns Hopkins University Applied Physics Laboratory \\
            Laurel, U.S.A\\
            \email{James.Spall@jhuapl.edu}
}

\date{Received: date / Accepted: date}

\maketitle

\begin{abstract}
Many public health threats exist, motivating the need to find optimal intervention strategies. Given the stochastic nature of the threats (e.g., the spread of pandemic influenza, the occurrence of drug overdoses, and the prevalence of alcohol-related threats), deterministic optimization approaches may be inappropriate. In this paper, we implement a stochastic optimization method to address aspects of the 2009 H1N1 and the COVID-19  pandemics, with the spread of disease modeled by the open source Monte Carlo simulations, FluTE and Covasim, respectively. Without testing every possible option, the objective of the optimization is to determine the best combination of intervention strategies so as to result in minimal economic loss to society. To reach our objective, this application-oriented paper uses the discrete simultaneous perturbation stochastic approximation method (DSPSA), a recursive simulation-based optimization algorithm, to update the input parameters in the disease simulation software so that the output iteratively approaches minimal economic loss. Assuming that the simulation models for the spread of disease (FluTE for H1N1 and Covasim for COVID-19 in our case) are accurate representations for the population being studied, the simulation-based strategy we present provides decision makers a powerful tool to mitigate potential human and economic losses from any epidemic. The basic approach is also applicable in other public health problems, such as opioid abuse and drunk driving.

\end{abstract}

\section{Introduction}
\label{sec1}
Numerous public health problems (e.g. global pandemic influenza, prescription drug overdoses, and obesity) continuously threaten people's lives, and finding optimal intervention strategies for many of these threats can be formulated as an optimization problem. As an example, by limiting the amount of opioids prescribed, regulating drug prices, and lecturing about the dangers of overdosing, professionals aim to find values that minimize costs related to opioid overdoses. Moreover, public safety entities can optimize speed limits on highways, locations of traffic lights, and traffic light cycles to find settings that minimize motor vehicle accidents. Solving public health optimization problems presents several challenges, including devising objective functions to compare the effectiveness of different strategies and utilizing suitable algorithms to optimize results under noisy measurement values since the inherent randomness involved. The main goal of this paper is to demonstrate some applications of a discrete stochastic optimization algorithm in determining optimal intervention strategies to yield minimal economic losses in the face of public health threats. Our focus below is on two infectious disease: influenza and COVID-19.\newline
\hspace*{10pt} Researchers in the field of epidemiology have attempted to solve optimization problems related to pandemics in several ways. Prosser et al. [23] evaluated the cost-effectiveness of the H1N1 vaccination across the United States in terms of the incremental cost-effectiveness ratio (ICER) and reported that vaccinating only children and working-age adults is more cost-effective compared to other vaccination strategies. Khazeni et al. [15] measured the cost-effectiveness of H1N1 vaccinations for New York City by infections, death averted, and qualified adjusted life-years, revealing that early vaccination is more cost effective relative to late vaccination. Piguillem and Shi [22] understood the optimal response to COVID-19 as a welfare maximization problem of a planner, and observed that extreme measures like mandatory quarantines seemed optimal for various curvatures of the welfare function. They also found that testing is a very close substitute of quarantine and can substantially reduce the need for indiscriminate quarantines. Hoertel et al. [13] evaluated the potential impact of quarantine duration, quarantine lifting type, post-quarantine screening, and use of a hypothetical effective treatment against COVID-19 with a microsimulation model. They projected that a two-step lifting of a quarantine according to age substantially lowered cumulative mortality and incidence. Charpentier et al. [7] conducted a tractable quantitative analysis of the impact of lockdowns and detection interventions using an extended SIR model (Susceptible, Infectious, and Recovered compartment model) for COVID-19, identifying the optimal lockdown policy to be an intervention structured in four successive phases. \newline
\hspace*{10pt} In the studies reviewed above, researchers mainly conducted sensitive analyses over a small range of candidate intervention strategies chosen on the basis of their prior knowledge and then selected the optimal one. Some researchers implemented a numerical algorithm to find the optimal set of intervention parameters, but that was based on a model with transmission dynamic of the pandemic governed by some known deterministic functions [24]. As discrete event simulation (DES) models can test assumptions at low cost and observe potential outcomes of decisions prior to implementation [12], they can be useful tool for our goal to develop improvement strategies. Among them, stochastic agent-based simulations can incorporate individual characteristics and person-to-person contact networks. Deterministic optimization algorithm will no longer be proper for such models as we only observe noisy measurement values and the dynamics can be complex and unknown. Paleshi et al. [20] implemented an agent-based model to simulate the spread of disease and optimized with a ranking and selection procedure called NSGS [19] to find a strategy on school closure and home confinement that minimized the effects of pandemic influenza. However, with the development of computational epidemiology, more factors can be taken into consideration in the construction of agent-based models. Simulating every feasible solution in ranking and selection algorithms becomes too computationally costly for high dimensional problems. Therefore, we employ a stochastic model and stochastic optimization method to estimate the optimal strategy among all possible options without testing every single one.\newline
\hspace*{10pt} Our goal to determine the optimal strategy motivates the question of how to measure the performance of an intervention strategy. We use a dynamic model built by [5] in our application concerning 2009 H1N1 and an agent-based model developed by [14] in our application concerning COVID-19. These two models simulate the spread of those two infectious diseases in a given population using transmission parameters estimated from historical data. The performance of each input intervention strategy is evaluated by its respective simulation outcome. The models are Monte Carlo simulations to reflect the uncertainties surrounding the transmission of infectious disease. Therefore, optimization of these kinds of problems can only be performed with noisy measurements of the actual performance. In addition, because the decision variables are discretely valued, the discrete simultaneous perturbation stochastic approximation method (DSPSA) is suitable for this problem. An early study conducted by [36] obtained the optimal intervention strategy by using DSPSA in simplified settings. In contrast, one modification of this report is the manner in which the intervention parameters are transformed in order to achieve better convergence in a relative small number of iterations. We also use an improved approach for computing the influenza-related cost and update the values into dollars of current purchasing power aligned with the Consumer Price Index [4]. As noted in [34], non-constrained cases are rare in stochastic optimization, so we address the constraints with appropriate projection functions. We will discuss the projection function that is used to address the parameter constraints in Section 2 of this paper.\newline
\hspace*{10pt} Our paper is organized as follows. First, we describe the DSPSA algorithm and discuss how to choose the gain sequence in a discrete high noise setting. Next, we introduce an application of the algorithm to the 2009 H1N1 pandemic, including the choice of simulator, the construction of the loss function, the transformation of input parameters, and some numerical results. Then we show a second application of the algorithm to the COVID-19 pandemic in the same manner. Finally, we discuss the potential real world interpretation of the final solution and offer ideas about how similar approaches can be applied in analyzing the optimal intervention strategy for other public health problems. While this paper is not oriented toward breaking new methodological ground, it does demonstrate some of the nontrivial considerations involved in the application of a powerful method in stochastic optimization in an important real-world problem.

\section{Constrained DSPSA Algorithm }
\label{sec2}
\hspace*{10pt} The DSPSA algorithm was introduced to solve discrete stochastic optimization problems in [35]. Suppose $\boldsymbol{\uptheta}\in\mathbb{Z}^{p}$, where $\mathbb{Z}=\{\cdots,-1,0,1,\cdots\}$. Our goal is to find a parameter $\boldsymbol{\uptheta}$ that minimizes the loss function $\textit{L}(\boldsymbol{\uptheta})$, given some constrained domain $\boldsymbol{\Uptheta}$; i.e., 
\begin{equation}
\min_{\boldsymbol{\uptheta}\in \boldsymbol{{\Uptheta}}}{\textit{L}(\boldsymbol{\uptheta})}
\end{equation}
Adding a restriction is equivalent to adding a projection function that maps the original $\boldsymbol{\uptheta}$ to a subset $\boldsymbol{\Uptheta}\subset\mathbb{Z}^{p}$. When the constraint is a well-defined range, we can easily define the projection $\boldsymbol{\Uppsi}$ as a piecewise linear function. For example, when each component of $\boldsymbol{\uptheta}=(\textit{t}_{1},\textit{t}_{2},\ldots,\textit{t}_{p})^{T}$ has a lower bound and an upper bound, $\textit{l}_{i}\leq\textit{t}_{i}\leq\textit{u}_{i},\textit{i}=1,\ldots,\textit{p}$, the projection can be represented as 
\begin{equation}
\boldsymbol{\Uppsi}(\textit{t}_{i}) =
\left\{
	\begin{array}{ll}
		\textit{l}_{i}  &\hspace{10pt}\textit{t}_{i}<\textit{l}_{i} \\
		\textit{t}_{i} &\hspace{10pt}\textit{l}_{i}\leq\textit{t}_{i}<\textit{u}_{i} \\
		\textit{u}_{i}-\uptau &\hspace{10pt}\textit{t}_{i}\geq\textit{u}_{i}
	\end{array}
\right.
\end{equation}
where $\uptau>0$ is a very small number, and is used to ensure the correctness of the DSPSA algorithm. Suppose $\textit{y}(\boldsymbol{\uptheta})$ is a noisy measurement of $\textit{L}(\boldsymbol{\uptheta})$; i.e., $\textit{y}(\boldsymbol{\uptheta})=\textit{L}(\boldsymbol{\uptheta})+\upepsilon(\boldsymbol{\uptheta})$, $\upepsilon(\boldsymbol{\uptheta})$ is the noise. The DSPSA method uses $\textit{y}(\boldsymbol{\uptheta})$ to find the minimum of $\textit{L}(\boldsymbol{\uptheta})$. Wang [34] gives a complete description of the theory and extends DSPSA to the case of a bounded domain. The algorithm is described below:\newline
\hspace*{10pt}Step 1: Pick an initial guess $\hat{\boldsymbol{\uptheta}}_{0}\in\mathbb{Z}^{p}$, and set $\textit{k}=0$.\newline
\hspace*{10pt}Step 2: Generate the random perturbation vector $\boldsymbol{\Updelta}_{k}=(\Updelta_{k1},\Updelta_{k2},\cdots,\Updelta_{kp})^T$, where the $\Updelta_{ki}$ are independent Bernoulli random variables taking values $\pm1$, each with probability 1/2, in this project. (Certain non-Bernoulli discrete distributions for $\boldsymbol{\Updelta}_{k}$ may also be used as long as $\hat{\boldsymbol{\uptheta}}_{k}^{\pm}$ in step 4 takes on valid integer values.)\newline
\hspace*{10pt}Step 3: Let $\boldsymbol{\uppi}(\hat{\boldsymbol{\uptheta}}_{k})=\lfloor\boldsymbol{\Uppsi}(\hat{\boldsymbol{\uptheta}}_{k})\rfloor+\textbf{1}_{p}/2$, where $\textbf{1}_{p}$ is a $\textit{p}$-dimensional vector with all components being unity, and $\lfloor\boldsymbol{\Uppsi}(\hat{\boldsymbol{\uptheta}}_{k})\rfloor$ denotes the floor of all components (round down to nearest integer). Note that $\boldsymbol{\uppi}(\hat{\boldsymbol{\uptheta}}_{k})$ is the middle point of a hypercube that contains $\hat{\boldsymbol{\uptheta}}_{k}$ with all corner points lying in $\mathbb{Z}^{p}$ \newline
\hspace*{10pt}Step 4: Evaluate $\textit{y}$ at the points $\hat{\boldsymbol{\uptheta}}_{k}^{\pm}=\boldsymbol{\uppi}(\hat{\boldsymbol{\uptheta}}_{k})\pm\boldsymbol{\Updelta}_{k}/2$, and form the estimate of the sub-gradient $\hat{\boldsymbol{g}}_{k}(\hat{\boldsymbol{\uptheta}}_{k})$ such that
\begin{equation}
\hat{\boldsymbol{g}}_{k}(\hat{\boldsymbol{\uptheta}}_{k})=(\textit{y}(\hat{\boldsymbol{\uptheta}}_{k}^{+})-\textit{y}(\hat{\boldsymbol{\uptheta}}_{k}^{-}))\boldsymbol{\Updelta}_{k}^{-1},
\end{equation}
where $\boldsymbol{\Updelta}_{k}^{-1}=(\Updelta_{k1}^{-1},\Updelta_{k2}^{-1},\cdots,\Updelta_{kp}^{-1})^{T}$ and $[\cdot]$ is the round operator applied to all components of the argument vector. \newline
\hspace*{10pt} Step 5: Update $\hat{\boldsymbol{\uptheta}}_{k+1}$ according to 
\begin{equation}
\hat{\boldsymbol{\uptheta}}_{k+1}=\hat{\boldsymbol{\uptheta}}_{k}-\textit{a}_{k}\hat{\boldsymbol{g}}_{k}(\hat{\boldsymbol{\uptheta}}_{k}) 
\end{equation}
where $\textit{a}_{k}=\textit{a}/(1+A+k)^{\upalpha}$ is the gain sequence, $\textit{a}>0, \textit{A}\geq0$, and $0.5<\upalpha\leq 1$. Either terminate the algorithm or set $\textit{k}=\textit{k}+1$ and return to Step 2. \newline
\hspace*{10pt} Step 6: After the terminal iteration (say, $k=M$), set the approximated optimal solution to be $[\boldsymbol{\Uppsi}(\hat{\boldsymbol{\uptheta}}_{M})]$. \newline
\hspace*{10pt} Suppose that $\boldsymbol{\uptheta}^{*}$ is the unique minimal point of $\textit{L}(\boldsymbol{\uptheta})$. By [35], under some common conditions of stochastic approximation algorithms, strong convergence of $\hat{\boldsymbol{\uptheta}}_{k}$ to $\boldsymbol{\uptheta}^{*}$ holds for the algorithm described above. Therefore, given a credible Monte Carlo simulation, we can obtain the optimal intervention strategy that minimizes the societal cost from the algorithm even though we only have noisy measurement values from a simulation.

\section{Application to 2009 H1N1}
\label{sec 3}
\subsection{Choice of epidemic model}
\label{sec 3.1}
In this application, we use a real-time agent-based influenza model (FluTE) developed by [5] to simulate the spread of influenza across major metropolitan areas. An agent-based model focuses on individual behavior. Individual agents with different demographic characteristics are programmed to behave and interact with other individuals in different ways, which makes the model more realistic than models based on aggregated responses. \newline
\hspace*{10pt}FluTE is freely available at \href{https://www.cs.unm.edu/~dlchao/flute/}{https://www.cs.unm.edu/$\sim$dlchao/flute/}. It provides high flexibility over a variety of factors, including pre-existing immunity, vaccine efficacy, response delay, vaccination fraction, and antiviral policy, among other factors. The relationship of these factors to one another is of such complexity that real-time simulation models offer a better approach to investigate the population-level effects of mitigation strategies in comparison to analytic models.\newline
\hspace*{10pt}In order to obtain the noisy measurement $\textit{y}$, we first run FluTE with the perturbations of current estimate $\hat{\boldsymbol{\uptheta}}_{k}$, which produces a summary in text format. The summary contains demographic features of the simulated population, a record of the chosen intervention strategies, and the outcome statistics, including the number of vaccines and antiviral agents used and the total number of symptomatic individuals corresponding to different age and risk groups. We then calculate monetary loss using these outcomes.

\subsection{Choice of Optimization Parameters and Transformation}
\label{sec 3.2}
We choose the optimization parameters to be vaccination fraction $(\mathit{F})$, vaccination priorities (a row vector, $\boldsymbol{\mathit{P}}$), antiviral policy $(\mathit{A})$ and school closure weeks $(\mathit{S})$. As DSPSA only works with integer vectors, we discretize and map the decision parameters $\boldsymbol{\uptheta}=(\mathit{F},\boldsymbol{\mathit{P}},\mathit{A},\mathit{S})^T$ into integer vectors. Vaccination fraction $(\mathit{F})$ is a scalar parameter taking real values between 0 and 1 and corresponding to the fraction of people to be vaccinated among each population group that is assigned vaccines. For economic efficiency, we test for values of this fraction less than 1 as we hope to achieve herd immunity by vaccinating only a fraction of population. We discretize the interval [0,1] into 0, 0.1, 0.2, $\ldots$ , 1.0, divide the values by 0.1 and convert them into integers 0, 1, 2, $\ldots$ , 10. This parameter will be mapped back to the feasible domain of [0,1] when we run FluTE. Vaccine priorities ($\boldsymbol{\mathit{P}}$) is a row vector of five integers, representing the vaccine priority for five categories of individuals. Though FluTE enables differential vaccine priorities setup for at most thirteen groups of individuals, some groups are included as part of other groups; for example, the pregnant women group is a subset of the young adults group (ages 19\textendash29) and the older adults group (ages 30\textendash64). Such overlap violates the basic assumption for the convergence of DSPSA, which requires $\textit{L}(\boldsymbol{\uptheta})$ to have a unique minimal point. Also for the goal of understanding which age group affects the efficiency of the intervention strategy the most, we only focus on the five non-overlapping categories that are grouped solely based on age, setting no specific priority for the other eight categories, 
which corresponds to setting their priorities as 0 by the software. To facilitate convergence, we have 3 indicate the highest priority, 2 be the next-highest priority, 1 indicates the third-highest priority, and 0 represent no vaccination for that group. For antiviral policy $(\mathit{A})$ and school closure weeks $(\mathit{S})$, we use the same rule as that of [35], where $\mathit{A}\in\{0,1,2,3\} $ (0: ``none", 1: ``treatment only", 2: ``HHTAP100", 3:``HHTAP", where HHTAP means that household members all get prescribed drugs if one member is ascertained and HHTAP100 to be a special option for LA county: drugs can go to the first 100 households that have a member ascertained), $\mathit{S}\in\{0,1,\cdots\}$ represents the number of weeks for school closure. 

\subsection{Loss Function}
\label{sec3.3}
Although vaccination, antiviral medicines and school closure are all effective for influenza prevention and mitigation, we must also consider the cost they incur. This guides us in the construction of loss function, which should contain the investment on prevention and treatment and health benefit loss incurred from infection as well as the possible hidden loss (e.g. cost from side effects of vaccines). We formulate the loss function in the following general format:
\begin{equation*}
\begin{split}
\textit{L}(\boldsymbol{\uptheta})&=E[\text{Medication Cost}(\boldsymbol{\uptheta})+\text{Vaccination Cost}(\boldsymbol{\uptheta})+\text{Antiviral Cost}(\boldsymbol{\uptheta})\\
&+\text{School Closure Cost}(\boldsymbol{\uptheta})+\text{Death Cost}(\boldsymbol{\uptheta})]
\end{split}
\end{equation*}
\hspace*{10pt}Note that at a fixed value of $\boldsymbol{\uptheta}$, the argument inside the $[\cdot]$ above is random due to the inherent unpredictability in cost response. The expectation operator forms the average over these random costs at the specified value of $\boldsymbol{\uptheta}$. Medication consists of two parts: non-hospitalized and hospitalized medications for influenza and its complications. Non-hospitalized medication cost is defined to be the summation of the cost of over-the counter (OTC) medications, physician visit cost for uncomplicated influenza (cost for prescription drugs included), physical visit cost for otitis media, and physician visit cost for non-hospitalized pneumonia, where otitis and pneumonia are complications of influenza. We assume both that all symptomatic individuals will use OTC medicine and that the age distribution in our synthetic city is proportional to the 2000 U.S. census [28], on which FluTE is based and which differs no more than $2.5\%$ from the population makeup in the 2010 census [29] for our age groups of interest.  The probabilities of the above four types of non-hospitalization medication and their corresponding costs of each age group are available in Tables S1 and Table S2 of [23]. These tables also provide different percentages of infection and costs for high risk and low risk people. We follow their guidelines accordingly, but some of their age group categories are different from that of the FluTE. Therefore, we transform the data provided in [23],  with age categories adjusted according to the age distribution provided by [28]. \newline
\hspace*{10pt} We calculate the hospitalized medication cost as the expectation of the total expense for all infected people who required hospital or ICU care. Chao et al. [6] suggests the division of each age group of symptomatic people into two risk-based subgroups: high risk and low risk, where symptomatic high risk individuals have a higher rate of being hospitalized. We employ the hospitalization ratio shown in Web Table 7 of [6] and follow [15]'s estimate that hospitalized individuals required five days of hospital care and that 10$\%$ of hospitalized patients required ten days of ICU care. According to Table 1 in [15], the average cost of hospitalization and ICU are $\$2430$ and $\$4960$ per day per person, respectively, when adjusted to 2019 dollars by multiplying the dollar amount by the ratio of 2019 CPI for medical care and that of the year when the referenced amount is based [4].\newline
\hspace*{10pt} Vaccination cost includes the cost of vaccines and the cost of vaccination-related adverse events. The cost of vaccines is estimated as $\$40$ (in 2019 dollars) per dose [36]. On the other hand, vaccines can also incur medication costs from vaccination-related adverse events, including systematic reactions, injection site reactions, anaphylaxis, and Guillain-Barr syndrome [23]. The probability and cost of each event are estimated in Tables S1 and Table S2 of [23]. The adjustment procedure is similar to what we did to the non-hospitalized medication cost.\newline
\hspace*{10pt} The total cost of antiviral agents is equal to the antiviral cost per dose times the number of doses used. The antiviral cost per dose includes the production cost and dispensing cost, which is approximately $\$74$ (in 2019 dollars) per course [11]. The simulator outputs the total number of antiviral agents used into the summary file. Multiplying this number by 74 gives the antiviral treatment cost.\newline
\hspace*{10pt} School closure cost is defined to include the cost of making up classes and the cost of the parents' lost wages to take care of their children. We employ the cost for school closure per day per students given by [11], which is $\$23$ (in 2019 dollars), as the cost for making up classes. For the cost incurred on parents' lost wages, Halder et al. [11] estimate the average wage of one person to be $\$980$ (in 2019 dollars) weekly and Araz et al. [1] estimate the number of days of work missed to be 2.5 for couples and 5 for single parents. Therefore, school closure cost per student per day is calculated using the following formula:
\begin{equation*}
\begin{split}
\text{School Closure Cost}&=\text{ Cost of Make Up Classes + Daily Wage of Parent} \\
&\times\text{ Days Work Missed }\times \frac{1}{5}
\end{split}
\end{equation*}
\hspace*{10pt} The final school closure cost per student per day is calculated to be $\$123$ (after adjusting to 2019 dollars). We then multiply this number by the number of students in each community and by 5 to get the weekly school closing cost per community. Chao et al. [5] specify the number of students in each community in FluTE. The number of communities can be retrieved from the simulation summary file. The resulting school closing cost per week per community is $\$221,804$, in 2019 dollars after we adjust the dollar amount by the ratio of 2019 General CPI to that of the year when the referenced amount is based [4]; therefore, we expect the optimal school closing days to be very close or equal to 0 weeks (0 days).\newline
\hspace*{10pt} We define the mortality loss to be the expected death cost, which has the following formula: 
\begin{equation*}
\begin{split}
\text{Mortality Loss }&=\text{ Number of Symptomatic Individuals}\\ &\times\text{ Fatality Ratio of Symptomatic Individual }\times\text{ Death Cost}
\end{split}    
\end{equation*}
The simulator outputs the total number of symptomatic individuals by age and the number of symptomatic high-risk individuals by age into the summary file; both are vectors. By subtracting the latter from the former, we can get the number of symptomatic non-high-risk individuals by age. Chao et al. [6] give the fatality ratio by age and risk group in Web Table 7. This table also shows that high-risk individuals have a higher fatality ratio than non-risk individuals.\newline
\hspace*{10pt} Molinari et al. [17] give the mortality loss in terms of both the value of a statistical life (VSL) and the present value of future earnings. The VSL estimates the number of dollars people are willing to pay for a certain number of reductions in mortality risks, whereas the present value of future earnings reflects how much output society loses in dollar terms when death occurs. In our context, the present value of future earnings better represents the definition of death costs as it reflects the economic loss to society. Therefore, we employ the data provided in Table 2 of [17], which are the present values of future earnings in 2003 dollars, and we adjust it to 2019 dollar values using historical CPI provided by [4].

\subsection{Numerical Results}
\label{sec3.4}
In this section, we present the numerical results of the performance of DSPSA algorithm. The general setup is as follows. Our synthetic population consists of 99,617 people from 20 contiguous tracts of central Los Angeles County. To facilitate the possible duplication of our experiment, we list the FIPS codes for the relevant tracts (FIPS codes specify the tract within Los Angeles County chosen for the simulation) in Appendix B. For vaccine availability, Chao et al. [6] gives a detailed U.S. pandemic H1N1 vaccine supply in its Web Table 3. We follow their assumption that Los Angeles County receives 3.195$\%$ of the total supply for the United States and that vaccine allocations after November 27, 2009 should be ignored. In addition, we assume that the vaccine supply for our synthetic city is proportional to the Los Angeles supply by population. We also set the initial inventory for all of the nine vaccines to be 0. Chao et al. [6] argues that all vaccines will be available in Los Angeles County 9 days after the U.S. supply date. Thus, we construct a detailed delayed daily vaccine supply table for our synthetic city, assuming that the simulation starts on September 1, 2009 and ends 175 days (25 weeks) later.\newline
\hspace*{10pt}  The initial value for the optimization process is that no intervention policy is needed. This value is expressed as $\hat{\boldsymbol{\uptheta}}_{0}=(2,0,0,0,0,0,0,0)^{T}$, which is equivalent to $(0.5,0,0,0,0,0,none,0)^{T}$ in terms of FluTE input. Following principles in section 7.5 of [25], section 3.3 of [34] provides the guidelines for coefficient selection in the gain sequence $\mathit{a}_{k}=\mathit{a}/(1+A+k)^{\upalpha}$. $\mathit{A}$ is recommended as 10$\%$ of the total number of iteration; $\upalpha$ is recommended to take a value approximately as small as formally allowed, and the value of $a$ can be found numerically by making the multiplication of $a$ and the average magnitude of $\hat{\boldsymbol{g}}_{0}(\hat{\boldsymbol{\uptheta}}_{0})$ equal to the desired magnitude of change of $\hat{\boldsymbol{\uptheta}}_{k}$ in early iterations. We set our number of iterations to be 10,000 and, following the guidelines above, we set $\mathit{A} = 1000$ and $\upalpha = 0.501$. For coefficient $\mathit{a}$, we find an initial estimate $\mathit{a}=1.5$.\newline
\hspace*{10pt} Before running the algorithm for a large number of iterations, we test the common random numbers (CRN) variance reduction technique [25, Chapter 14] for a small number of iterations to determine whether it can improve the performance. Plots of trials reveal no significant improvement when using the CRNs. The lack of effectiveness aligns with separate numerical experiments showing that no statistically significant positive correlation exists between $\textit{y}(\hat{\boldsymbol{\uptheta}}_{k}^{+})$ and $\textit{y}(\hat{\boldsymbol{\uptheta}}_{k}^{-})$ when using the same random number seed, a result that renders the CRNs ineffective. A lack of synchronization in the two outputs $\textit{y}(\hat{\boldsymbol{\uptheta}}_{k}^{+})$, $\textit{y}(\hat{\boldsymbol{\uptheta}}_{k}^{-})$ from FluTE may be a possible reason for the lack of correlation. To bring this to light, after some number of iterations, $\hat{\boldsymbol{\uptheta}}_{k}^{+}$ may indicate that only one age group is prioritized in receiving vaccines, while $\hat{\boldsymbol{\uptheta}}_{k}^{-}$ indicates that two groups are prioritized. Then by the design of the model, we need different numbers of random variables in their simulations and fixing the same seed can only ensure partial alignment. Besides, we do not have the contribution of the perturbation vector $\boldsymbol{\Updelta}_{k}$ converging to zero in DSPSA, which is one of the conditions that guarantees that CRNs aid in continuous variable versions of the SPSA algorithm. As a result of little synchronization and the perturbation size not converging to zero, we do not incorporate the CRNs technique into DSPSA in this application.\newline
\hspace*{10pt} We first check the reasonableness of applying the DSPSA algorithm by examining the trends in total loss over the iterations. Fig. 1 illustrates the tendency of total loss over 10000 iteration steps averaged over 10 Monte Carlo trials. As there is an approximately $60\%$ decrease in total loss relative to a strategy of no intervention, the DSPSA algorithm is demonstrated to be suitable for seeking an optimal resource allocation in the epidemiological problem of interest here.\newline
\hspace*{10pt} We then check the sensitivity of our simulation. In each trial, the solution greatly reduced the total monetary cost from the initial guess. However, as the algorithm is based on noisy measurement of $\textit{L}(\boldsymbol{\uptheta})$, the noise level largely determines the reliability of the result. Fig. 2 shows the fluctuation of a single trial, from which we can see that precise interpretation of the result is a problem with a high-noise setting. If one strategy has a loss value $L$ only slightly below that of another strategy, such superiority might not be detected by the algorithm in limited iterations. As a result, we may only obtain a set of good choice strategies, but it is difficult to discern the best strategy.

\begin{figure}
\centering
\includegraphics[scale=0.4]{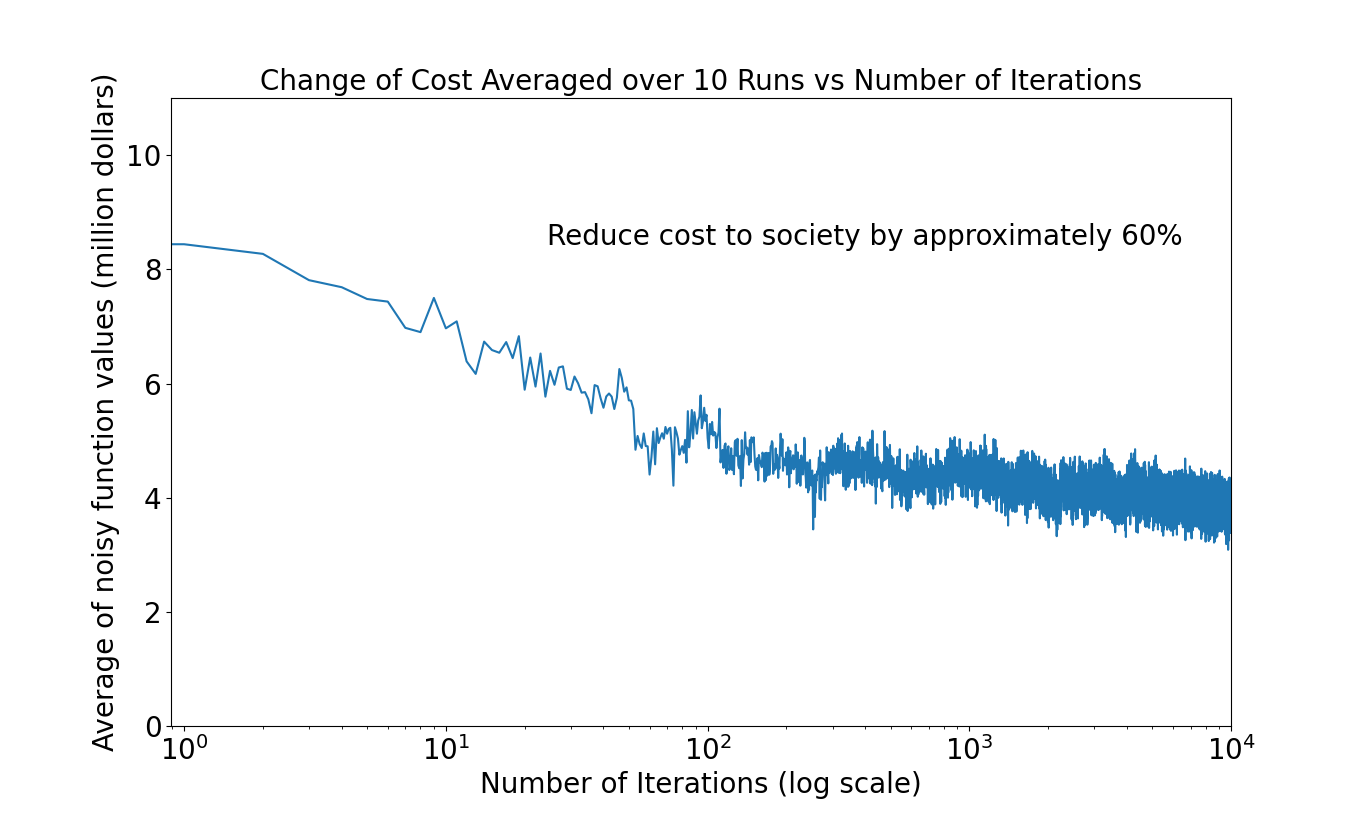}
\caption{Change of Cost Averaged over 10 Runs Versus Number of Iterations}
\end{figure}

\begin{figure}
\centering
\includegraphics[scale=0.4]{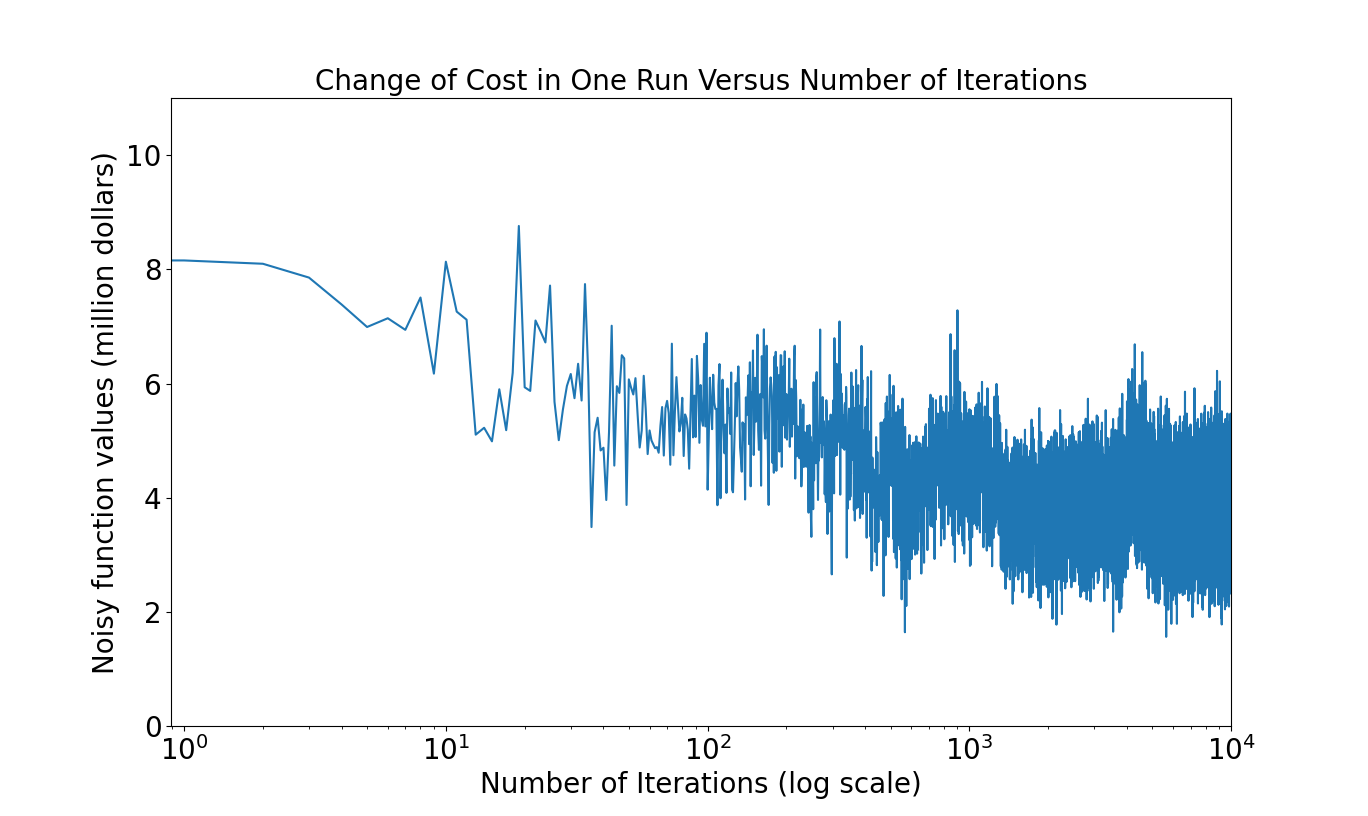}
\caption{Change of Cost in One Run Versus Number of Iterations}
\end{figure}

\indent From the trials, we find that the optimal vaccination strategy is to vaccinate solely the school-age children and to vaccinate all of them. This result agrees with the conclusion of [3], where vaccinating school-age children is the most effective way to reduce the overall influenza attack rate as well as the overall number of deaths. The result is derived as follows, given the condition that pre-school, young adults, older adults, and the elderly will receive no vaccine and that only school age children will receive vaccines, then we should target to vaccinate $100\%$ of them. Our purpose of setting a vaccination fraction is to vaccinate only a proportion of the population to achieve herd immunity. Therefore, we want to verify whether this proportion number follows the basic herd immunity threshold formula. From [6], we know that the basic reproductive number $\textit{R}_{0}$ of H1N1 is 1.3, then the herd immunity threshold for total population can be calculated as: $$1-\frac{1}{\textit{R}_{0}}=1-\frac{1}{1.3}=23.1\%$$ 
This number means that $23.1\%$ of the total population should be vaccinated to achieve herd immunity. Comparing to the optimal value above, the proportion of population that actually gets vaccinated in the simulation of Fig.2 is: $$\frac{\text{number of school age children}}{\text{total number of population}}=\frac{21,976}{99,617}=22.1\%$$
\hspace*{10pt} The simulation optimization result is close to the analytical solution. We see that this is the particular case when the school-age population is the dominating factor of H1N1 influenza transmission in a small population setup, and the optimal proportion of population to vaccinate is that entire age group. This agrees with the argument in [16] that pushing high vaccination coverage on school-age children is essential to building herd immunity for H1N1 influenza. Among four possible choices of antiviral policy, the optimal solution should be HHTAP. This result is consistent with the result in [36], where they also assume an unlimited antiviral agent supply. For school closure days, our initial guess of 0 weeks (0 days) is the optimal solution as closing school cost is significant compared to other parameters.\newline
\hspace*{10pt} To determine the robustness of our solution, we carry out 500 independent Monte Carlo simulation trials to get the noisy loss function values at the terminal iterate. Then we estimate the mean and the variance of the terminal loss function value with the sample mean and sample variance, and build an estimated $95\%$ confidence interval for the mean loss function value from it, which is $[3.75\times10^{6}, 3.87 \times10^{6}]$. Across the 500 trials, the minimum noisy loss value is about 2.16 million, and the maximum noisy loss value is about 5.73 million. Then we compare the sample mean with average loss function values of optimal intervention strategies suggested by other public health papers applied to the same L.A. tracts. Halder et al. [11] supposed that a combination of antiviral policy (HHTAP) and two weeks school closure was effective for mitigating a pandemic influenza with 2009 H1N1 characteristics, and we compute the average noisy loss function value with such strategy from 500 simulations, which is larger than $\$$30.00 million. Prosser et al. [23] suggested that vaccinating school age children and adults of age 30\textendash40 would be most cost-effective, and the average loss function value for that strategy is about $\$$4.41 million. Pasquini-Descomps et al.[21] summarized that at least $60\%$ of the population should be vaccinated, and we get the average measurement value to be about $\$$5.76 million. Khazeni et al. [15] claimed that vaccinating $40\%$ of population would be cost-saving; then we compute that the average cost for such is about $\$$5.14 million. Though strategies suggested by other papers are not based on the same assumptions, we can see through the comparison that the optimal solution from DSPSA has significantly lower cost, which indicate that it is a reasonably good solution in the setting of this paper. This comparison provides additional confidence about the quality of our optimal solution.

\section{Application to COVID-19}
\subsection{Choice of epidemic model}
\label{sec 4.1}
We now consider the application of the above approach to COVID-19 prior to the availability of vaccines. We use an agent-based simulator (Covasim) developed by [14] to model the spread of the COVID-19 pandemic throughout a population. The model is freely available at \href{https://www.github.com/InstituteforDiseaseModeling/covasim}{https://www.github.com/InstituteforDiseaseModeling/covasim}, and it has some advantages over other available models for the COVID-19 pandemic. Covasim incorporates randomness in its simulation output, so it better captures the stochastic nature of the virus spread compared to popular deterministic epidemic models like the Susceptible-Infectious-Recovered (SIR) model. Covasim includes demographic information on age structure and population size, and it builds a transmission network with multiple social layers, including households, schools, workplaces, and communities to reflect individual contact in reality. Covasim also supports an extensive set of interventions, which works for our goal to explore scenarios under different combination of intervention strategies (pre-vaccine).\newline
\hspace*{10pt} To obtain the noisy measurement $y$, we first run Covasim with the current estimate $\hat{\boldsymbol{\uptheta}}_{k}$. Next, we collect some summary statistics from the simulation output, including the total number of diagnostic testing performed, cumulative number of symptomatic individuals, and others. We will then use these statistics to calculate the economic cost of the scenario, with all cost values listed in 2020 U.S. dollars.

\subsection{Choice of Optimization Parameter and Transformation}
\label{sec 4.2}
There are recent developments that would likely lead to future studies related to COVID-19. First, vaccines for the COVID-19 have been developed and are serving as powerful tools against the pandemic in many countries. Second, variants of the coronavirus (delta, omicron, etc.) have also spread throughout the world. Additionally, general strategy related to masking, social distancing, quarantining, and so on have been refined over time as public health authorities learn more about the disease. These changes affect our goal of making a factual representation of the current situation in the model. Therefore, for the evaluation of the simulation optimization method, we restrict our attention to the pre-vaccine scenario and do not include vaccination as an intervention strategy. In that sense, one may consider this a ``proof of concept" study, with vaccines and other interventions left for future study.\newline 
\hspace*{10pt} In this project, we choose optimization parameters that relate to four types of intervention strategies: social distancing ($\boldsymbol{D}$, including closure of workplace, cancellation of public events, etc.), school closure ($\boldsymbol{S}$), testing ($\boldsymbol{T}$), and contact tracing ($\boldsymbol{C}$). For each type of intervention, there are three parameters: the start and end day of the intervention and percentage level of intensity of the policy. Because DSPSA works only with integer vectors, we discretize the input vector $\boldsymbol{\uptheta}=(\boldsymbol{D},\boldsymbol{S},\boldsymbol{T},\boldsymbol{C})^{T}$ and map each parameter into integer values (each of $\boldsymbol{D},\boldsymbol{S},\boldsymbol{T},\boldsymbol{C}$ are row vectors). Note that dim$(\boldsymbol{\uptheta})=12$.  The percentage level of intensity for each intervention strategy takes a real value between 0 and 100. We discretize the level [0,100] into 0, 10, 20, $\cdots$, 100, divide the values by 10, and convert them into integers 0,1,2,$\cdots$, 10. These parameters will be mapped back to percentage levels in [0,100] when we run Covasim. The start and end day for each intervention strategy should be kept within the range of days of the simulation, so we make their lower bounds $l_{i}$ in the projection function (2) to be 1 and their upper bounds $u_{i}$ to be the length of days of the simulation. Under this setup, the number of possible interventions in $\boldsymbol{\Uptheta}$ would be about $10^{17}$. Taking into consideration the fact that the start day of one intervention should not be later than its end day, we add one more step following the projection function (2) in the algorithm: we set the start day to be one day before the end day whenever the end day of any policy turns out to be earlier than its start day after the projection function (2) in the algorithm.

\subsection{Loss Function}
\label {sec 4.3}
All four types of intervention strategies mentioned above are effective in helping control the spread of COVID-19. To determine the optimal combination of strategies, we construct a function to compute the socio-economic cost of each scenario. We formulate the loss function as follows:
\begin{equation*}
\begin{split}
\textit{L}(\boldsymbol{\uptheta})&=E[\text{social distancing cost}(\boldsymbol{\uptheta})+\text{school closure cost}(\boldsymbol{\uptheta})+\text{testing cost}(\boldsymbol{\uptheta})\\
&+\text{contact tracing cost}(\boldsymbol{\uptheta})+\text{treatment cost}(\boldsymbol{\uptheta})+\text{death cost}(\boldsymbol{\uptheta})]
\end{split}
\end{equation*}
\hspace*{10pt} All of the above costs are cumulative over the time period of interest. For the cost of social distancing, Strong and Welburn [26] estimate the percentage decrease in weekly household income for each state in the U.S. under several social distancing scenarios. Below we consider scenario 5B, which includes closing schools, bars, and restaurants, banning large events, closing nonessential businesses, and quarantining the most
vulnerable. Thunstrom et al. [27] assume that a set of social distancing measures similar to scenario 5B leads to an average 38$\%$ decrease in individual contact rate, a measure to reflect the transmission probability between an infected individual and a susceptible individual. Using the estimates of the reduction in contact rate and decrease in household income, we assume a linear relationship between the percentage level of intensity for social distancing and percentage decrease in weekly household income. Next, we collect census data about the population and number of households [30] to estimate the number of households for our synthetic population in simulation. Then we use the average weekly household income from the census data [30] to compute the social distancing cost by formula: 
\begin{equation*}
\begin{split}
\text{Social Distancing Cost}=&\frac{\text{Percentage level of intensity}}{38}  \times \text{Total household income} \\
&\times \text{Percentage decrease in weekly household income} \\
&\times \text{Duration of the policy in weeks.}
\end{split}
\end{equation*}
\hspace*{10pt} The total school closure cost is computed to be the sum of the cost of making up classes and the cost of
the lost wages for parents to take care of their children. We estimate the cost for making up classes with the school closure cost per day per student given by [11], which is $\$23$, in 2020 U.S. dollars. For the cost incurred on parents’ lost wages, Halder et al. [11] estimate the
average wage of one person to be $\$992$ (in 2020 dollars) weekly and Araz et al. [1] estimate the number of
days of work missed to be 2.5 for couples and 5 for single parents. Using these estimates, the total school closure cost per student per day is calculated to be $\$125$ (after adjusting
to 2020 dollars). Next, we use data from [30] to estimate the total number of students in our synthetic population. Using the closure cost per student and the total number of students, we get the daily school closure cost in our population under a complete school shutdown. Then the cost of the intervention policy is computed by multiplying the respective percentage level of intensity. \newline
\hspace*{10pt} Though COVID-19 tests are currently free for individuals in the U.S., they will incur a cost to health insurance companies or the government. We employ the cost of a test performed by the CDC as given by [32], which is $\$36$. Then we compute the total cost for testing by multiplying it by the total number of tests performed, which is collected from the simulation output.\newline
\hspace*{10pt} As for the cost of contact tracing, Watson et al. [37] claim that 100,000 tracers would be needed nationwide at an annual cost of $\$$3.6 billion. We then compute the number of tracers needed in our synthetic population and their sum cost by proportion. The total cost will finally be multiplied by the percentage level of intensity to obtain the cost of contact tracing policy.\newline
\hspace*{10pt} For the treatment cost, Bartsch et al. [2] estimate the medical costs for non-hospitalized patients with COVID-19 symptomatic infections to be $\$3,994$ and Cohen, Whittal, and Murray [9] estimate the average cost to treat a hospitalized patient with coronavirus to be about $\$30,000$. We assume that people who are symptomatic but have no severe or critical symptoms require no hospitalization, while those who have severe symptoms or critical symptoms, as well as those who die, all require hospitalization. Then we follow the previous estimates and compute their treatment costs correspondingly.\newline
\hspace*{10pt} To compute the cost of death, we use the Value of Statistical Life (VSL). Eichenbaum et al. [10] estimate that the average VSL for a person in the U.S. to be $\$$9.3 million, and we multiply their estimate by the total number of deaths from the simulation output to obtain the sum of death costs from the COVID-19 pandemic.

\subsection{Numerical Results}
\label{sec 4.4}
In this section, we present some numerical results relating to the application of the DSPSA algorithm. The general setup for our simulation is as follows. Our synthetic population consists of 100,000 people with the same age structure as the state of Maryland, as revealed by the census data [30]. We use the hybrid contact network in Covasim for our simulation, which will be suitable to represent the realistic contact network of Maryland people in general [14]. For parameters about the disease progression and transmission such as the recovery time for severe cases, we use the default values in Covasim as denoted in Section 2.2 and 2.3 of [14]. We set the start day of the simulation to be March 1, 2020 and the length of the simulation to be 60 days, with the initial number of infected people being 50. As a measure of lack of sensitivity to the initial number, the length of time for social distancing in the results change by less than 3 days on average when the number of initial infected people varies between 5 and 500.\newline
\hspace*{10pt} The initial value for the optimization process is that no intervention policy is implemented. This value is expressed as $\hat{\boldsymbol{\uptheta}}_{0}=(1,2,0,1,2,0,1,2,0,1,2,0)^{T}$. We follow the guidelines provided in Section 3.3 of [34] and set the gain sequence to be $a_{k}=a/(1+A+k)^{\upalpha}$, where $\textit{A}$ is 10$\%$ of the total number of iterations and $\upalpha$ is 0.501. The value of $a$ can be found numerically by making the multiplication of $a_{0}$ and the average magnitude of $\hat{\boldsymbol{g}}_{0}(\hat{\boldsymbol{\uptheta}}_{0})$ equal to the desired magnitude of change of $\hat{\boldsymbol{\uptheta}}_{k}$ in early iterations. We set our number of iterations to be 5,000 and follow the guidelines to find an initial estimate $\textit{a}=0.08$.\newline
\hspace*{10pt} Before running the algorithm for a large number of iterations, we perform some tests on the common random numbers (CRN) variance reduction technique [25, Chapter 14]. Numerical experiments indicate the existence of statistically significant positive correlation between $y(\hat{\theta}^{+}_{k})$ and $y(\hat{\theta}^{-}_{k})$ when using the same random number seeds, one condition that renders the effectiveness of the CRNs. Plots of trials also reveal improvement when using CRNs. As a result, we incorporate CRNs into the algorithm in this application.\newline
\hspace*{10pt} We first determine the reasonableness of applying the DSPSA algorithm in finding the optimal solution to this simulation by examining the fluctuation of total noisy loss $y(\boldsymbol{\uptheta})$ over the iterations (the true loss $L(\boldsymbol{\uptheta})$ is not available). Fig. 3 shows the tendency of total loss over the first 2,500 steps in a single trial. As there is approximately a 90$\%$ decrease in total noisy loss, DSPSA appears to be suitable for seeking optimal intervention strategy for COVID-19 in this simulation.\newline
\begin{figure}
\centering
\includegraphics[scale=0.4]{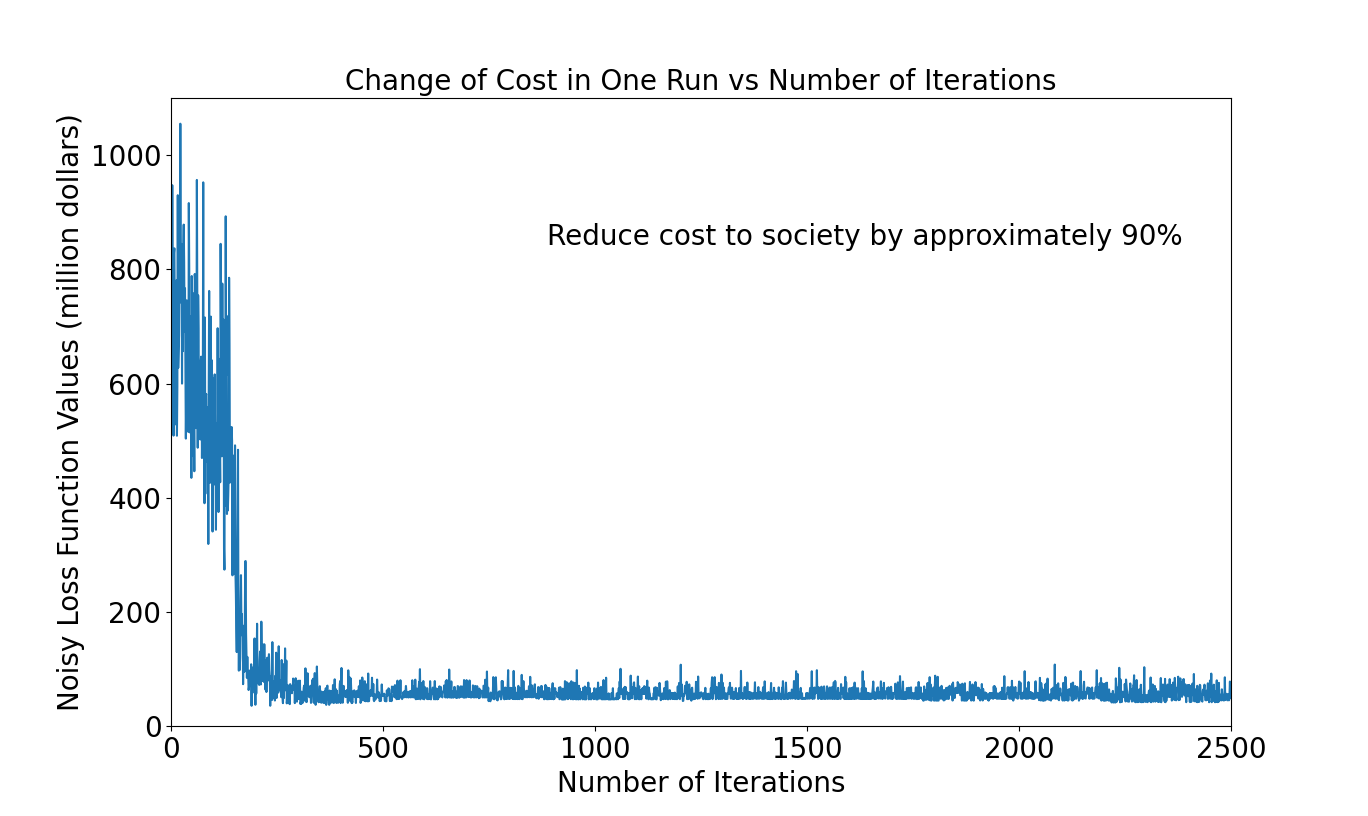}
\caption{Change of Cost in One Run Versus Number of Iterations}
\end{figure}
\hspace*{10pt} Next we check the sensitivity of our simulation results to alternative solutions. In each trial, the solution greatly reduced the total monetary cost from our initial guess. While the algorithm is based on noisy measurements of $L(\boldsymbol{\uptheta})$, the noise level largely determines the reliability of the result. From Fig. 3 we can see that the level of fluctuation decreases through the iterations. Therefore, we have confidence that the solution we get from the algorithm is a reasonable solution. However, due to the scale of the plot, there is a reasonable chance that one strategy might be superior to another by a tiny margin and such superiority may not be detected by the algorithm in limited iterations. It is, therefore, possible that the solution we find is only one from a set of good choice strategies, but perhaps not the best one. From Fig. 3 we also notice that after some number of iterations, the magnitude of fluctuation upward seems to be greater than that downward. This asymmetry appears to be related to the transmission dynamics of the disease. When the input intervention strategies are close to some optimal ones after some iterations, slight improvement in the intervention strategies is no longer able to significantly reduce the total number of infected people. However, slight deterioration of the strategies may cause asymmetrically more people to get infected by the end of the simulation through the high transmission rate of the SARS-CoV-2. Such increase in the number of symptomatic individuals increases the societal cost of the disease substantially, which determines the pattern in Fig. 3.\newline
\hspace*{10pt} From the trials, we find that the terminal solution is $\hat{\boldsymbol{\uptheta}}_{5000}=[4,18,10,1,2,0,1,20,10,1,2,0]^{T}$. The solution is interpreted to impose a social distancing policy with an extremely aggressive level of intensity (percentage level 100) for about two weeks starting from Day 4 of the simulation and to also implement massive testing of the general population (goal to test 100$\%$ of the population) starting from the first day for about 20 days. Then tests are performed to individuals who are symptomatic for the rest of the 2 months, and those who are tested positive are required to have two-week-home-based-quarantine. This solution relies on the fact that the length of time after exposure before an individual is infectious is set to have a mean of 4.6 days in Covasim [14], and it corresponds to suggestions from the CDC that the timeframe for self-quarantine of suspected COVID-19 cases should be about 14 days [31]. The recommendation of massive testing is also consistent with the results from [22] that testing generates sizable welfare gains. As for the suggested level of intensity for the social distancing policy, though it may not be possible to reach the extreme (100$\%$) level of intensity, and the number itself does not provide concrete ideas about specific measures that the government should impose, we can seek some guidance from the effect of established policies of other countries and nations. According to [38], the average number of contacts per individual in Wuhan, China decreased by about 86.3$\%$ through the duration of the social distancing policy imposed by the Chinese government, which corresponds to a social distancing policy of about 86 percent level of intensity by our definition. As their levels of intensity are close, we can learn that the desired social distancing policy suggested by the optimal solution should be close to that of Wuhan, China. Therefore, the government should impose similar measures including, but not limited to, compulsory wearing of masks in public areas, banning large events, and closing all nonessential business, among other measures [8].

\section{Conclusions}
\label{sec 5}
In this paper, we considered two infectious-disease-related problems with the goal of finding optimal mitigation strategies using the DSPSA method of multi-dimensional stochastic optimization for discrete problems. Numerical results reveal that the algorithm results in a significant decrease in loss after a reasonable number of iterations, which illustrates that the DSPSA algorithm is a good approach to solving optimization problems with noisy loss.\newline
\hspace*{10pt} Many previous public health studies analyze only a small set of candidate intervention strategies; such approaches require prior knowledge on the likely solution and are conducted on a case-by-case basis. Applications in this paper present a new approach that enables researchers to optimize over the whole set of available intervention strategies without the conceptual idea of a possible optimal strategy. This paper also sheds light on public health optimization analysis by considering not only the goal of public health but also the aim of minimal societal cost. Steps presented in this paper such as the loss function setup and the transformation of parameters can provide guidance in analyzing other types of public health problems, such as drug abuse (e.g., opioid) or further COVID developments (vaccines, antivirals, etc.) or variants of COVID (delta, omicron, etc.).\newline  
\hspace*{10pt} We may also conduct the same type of study on a modification of the DSPSA algorithm called the mixed simultaneous perturbation stochastic approximation (MSPSA), which applies to a mixture of discrete and continuous input parameters [33]. Because MSPSA is similar to DSPSA, only slight modifications are needed for applications from this paper, such as allowing the vaccine fraction ($\mathit{F}$) to range from [0,1] continuously. By the nature of the two algorithms, we expect that the solution will converge to an optimal solution in a similar manner as applications in this paper, and we leave the implementation for future study. As MSPSA formally allows for both continuous and discrete parameters, it can adapt to a broader range of public health problems.\newline
\hspace*{10pt}In summary, this paper is in response to a need to enhance response to public health threats at regional or larger scale. While current models can support decision makers by measuring the expected impact of factors such as school closures, vaccination campaigns, vaccine fractions, and antiviral treatment strategies, they do not automatically provide an optimal combination of the variables to achieve cost-effective societal aims. This paper demonstrated that a formal stochastic optimization method, DSPSA, combined with a user-selected Monte Carlo simulation of disease transmission and economic impact can provide solutions that are both practically realistic and potentially better than current ad hoc solutions. Results presented in this paper can also motivate us to apply the SPSA family (DSPSA or MSPSA) of algorithms in future public health mitigation studies such as developing countermeasure for bioterrorism. In the long term, the DSPSA/MSPSA-based approach can also adapt to improved simulation models, real-time data sources, better biological understanding, and new candidate intervention strategies (e.g., the current availability of vaccines for COVID-19) to help local and national authorities make better decisions.

\section*{Appendix}
Since all the 20 tracts we picked are located within California State and Los Angeles County, the State FPIS code and county FIPS code are 06 and 037 respectively for all tracts. The tract FIPS code for the selected 20 tracts are: 101110, 101120, 101210, 101220, 101300, 101400, 102101, 102102, 103101, 103102, 103200, 103300, 103400, 104103, 104104, 104105, 104106, 104107, 104201, 104202.

%
%

\begin{acknowledgements}
We thank Yan Zhou and Mengdan Zhang (the Johns Hopkins University) for all the preliminary work done during the course of this research.
\end{acknowledgements}

%
%



\end{document}